\documentclass[reqno]{amsart}[11pt]
\usepackage{amsmath, amssymb, latexsym, amsthm, amsfonts, mathrsfs, amscd, mathtools, commath, amscd, mathrsfs}
\usepackage{float}

\usepackage{hyperref}
\renewcommand{\eqref}[1]{(\eqref{#1})}   

\numberwithin{equation}{section}

\theoremstyle{plain}
\newtheorem{thm}{Theorem}
\newtheorem{lem}{Lemma}

\theoremstyle{definition}

\newcommand{\be} {\begin{equation}}
\newcommand{\ee} {\end{equation}}
\newcommand{\bea} {\begin{eqnarray}}
\newcommand{\eea} {\end{eqnarray}}
\newcommand{\ba} {\begin{array}}
	\newcommand{\ea} {\end{array}}
\newcommand{\nn} {\nonumber}

\begin{document}

\title[Zeros of some level N Eisenstein series]{A simple method to extract zeros of certain Eisenstein series of small level}

\author{Aradhita Chattopadhyaya}
\address{Department of Centre for High Energy Physics \\ Indian Institute of Science, Bangalore- 560012, India}
\email{aradhitac@iisc.ac.in} 


\begin{abstract}
{This paper provides a simple method to extract the zeros of some weight two Eisenstein series of level $N$ where $N=2,3,5$ and $7$. The method is based on the observation that
these Eisenstein series
 are integral over the graded algebra of modular forms on $SL(2,Z)$
 and their zeros  are `controlled' by those of $E_4$ and $E_6$ in the fundamental domain of $\Gamma_0(N)$.}
\end{abstract}
\maketitle 

\section{{\bf Introduction}}
In the work of Rankin and Swinnerton-Dyer \cite{RSD} the location of zeros of all Eisenstein series $E_k$ (of weight $k \geq 4$, even) of full modular group $SL_2(\mathbb{Z})$ had been determined. In the fundamental domain this was found to be always on the arc $|\tau|=1,\; {\rm with}\; 2\pi/3 \geq{\rm arg}(\tau)\geq \pi/2$. The method has been generalized to Fricke groups in recent works of \cite{shi07} and to the subgroups of $SL_2(\mathbb{Z})$ in \cite{gslt,gar1}. The zeros of  weight two Eisenstein series $E_2(q)$ were studied by \cite{woody}, \cite{bs}. In this work we shall find the zeros of Eisenstein series $\tilde{E}_N${\footnote{In general $\tilde{E}_N$ is defined as negative of what it is defined here in (\ref{en})}} which are holomorphic modular forms of weight 2 of $\Gamma_0(N)$ defined as
\be \label{en}
\tilde{ E}_N(\tau):=\frac{1}{N-1}(NE_2(N\tau)-E_2(\tau)),
\ee
\noindent
where $E_2(\tau)$ is the quasimodular Eisenstein series defined by,
\be
E_2(\tau):=1-24\sum_{n=1}^{\infty}\sigma_{1}(n) q^n,\qquad q=e^{2\pi i\tau},
\ee
with $\sigma_{1}(n)$ being the sum over all the divisors of $n$.

 The method which we present is quite different from that of \cite{RSD}, however can only be applied for $N=2,3,5,7$. 
\noindent
Our main observation is that $\tilde{E}_N$ is integral over the graded algebra
\[
 M(SL_2(\mathbb{Z}))=\bigoplus_{k\ge 0}M_k(SL_2(\mathbb{Z}))=\mathbb{C}[E_4,E_6].
 \] with $E_4,\, E_6$ being the Eisenstein series of weight 4 and 6 defined as,
 \begin{eqnarray} \label{e4}
 E_4(\tau):=1+240\sum_{n=1}^{\infty}\sigma_{3}(n) q^n, \\ \label{e6}
 E_6(\tau):=1-504\sum_{n=1}^{\infty}\sigma_{5}(n) q^n,
 \end{eqnarray}
 where $q=e^{2\pi i \tau}, \; \tau \in \mathbb{H}$.
 
 Let us denote the fundamental domain for $\Gamma_0(N)$ by $F_N$. For $N=2,3,5,7$ the zeros of $\tilde{E}_N$ are `controlled' (\ref{poly}) by those of $E_4$ and $E_6$ in $F_N$.

\noindent
We now state the main result of our paper which would be proved in section \ref{proof} .
\begin{thm}\label{theorem} 
All the zeros of $\tilde{ E}_N(\tau)$ in the fundamental domain of $\Gamma_0(N)$ lie
as follows:
		 \begin{table}[htb] \label{table}
		 	\renewcommand{\arraystretch}{0.5}
		 	\begin{center}
		 		\vspace{0.2cm}
		 		\begin{tabular}{|c|c|c|c|c|}
		 			\hline
		 			& & & &  \\
		 			$N$ & $2$ & $3$ & $5$ & $7$ \\
		 			& & & &  \\
		 			\hline
		 			& & & & \\
		 			$\tau$ &$-\frac{1}{i+1}$ & $-\frac{1}{e^{2\pi i/3}+2}$ & $-\frac{1}{i+2}, $& $-\frac{1}{e^{2\pi i/3}+3},$ \\
		 			& & & $-\frac{1}{i+3}$ & $-\frac{1}{e^{2\pi i/3}+5}$\\
		 			\hline
		 		\end{tabular}
		 	\end{center}
		 	\vspace{0.5cm}
		 	\caption{Zeros of $\tilde{ E}_N(\tau)$ for $N=2,3,5,7$.}
		 	\renewcommand{\arraystretch}{0.2}
		 \end{table}
	\end{thm}

Before we move into the proof of the theorem we shall first provide a brief description where these modular forms occur in the context of string theory in physics.

\section*{Physics Motivation}
Partition functions defined in string theory are generally related to modular functions. One such partition function defined over a $K3$ manifold is elliptic genus.
\begin{eqnarray} \label{ellip}
F(\tau, z) &=& {\rm Tr}_{RR} \left[ (-1)^{F_{K3} + \bar F_{K3} } e^{ 2\pi i z F_{\rm K3} } 
q^{L_0 - \frac{c}{24} } \bar q ^{\bar L_0 - \frac{\bar c}{24}} \right],  \\ \nonumber
&=& \sum_{b=0}^1 \sum_{j \in 2\mathbb{Z} + b, \; n \in \mathbb{Z}} 
c_b( 4n - j^2) e^{2\pi i n \tau + 2\pi i j z} .
\end{eqnarray}
The trace in the above equation is taken over the Ramond-Ramond sector of the 
${\mathcal N}=(4, 4)$ super conformal field theory of $K3$ with central charge 
$(6, 6)$ and $F$ refers to the Fermion number, $L_0$ and $\bar L_0$ are the scaling operators in the left moving and right moving CFTs.
The result of the above trace can be given by a weak Jacobi form of index 1 and weight 0 \cite{EOTY}.
\be
F(\tau,z)=8A(\tau,z)=8\left(\frac{\theta_2^2(\tau,z)}{\theta_2^2(\tau,0)}+\frac{\theta_3^2(\tau,z)}{\theta_3^2(\tau,0)}+\frac{\theta_4^2(\tau,z)}{\theta_4^2(\tau,0)}\right).
\ee
Here $\theta_i(\tau,z)$ are Jacobi theta functions.

 There are symplectic automorphisms of $K3$ related to the conjugacy classes of Mathieu group $M_{24}$ and we can define a partition function similar to the elliptic genus for these orbifolds of $K3$ called "twisted elliptic genus" of $K3$.
This is defined as
\begin{eqnarray} \label{expelipg}
F^{(r, s)}(\tau, z) &=& \frac{1}{N} {\rm Tr}_{RR\; g^{\prime r } }
\left[ (-1)^{F_{K3} + \bar F_{K3} } 
g^{\prime s} e^{ 2\pi i z F_{\rm K3} } 
q^{L_0 - \frac{c}{24} } \bar q ^{\bar L_0 - \frac{\bar c}{24}} \right],  \\ \nonumber
&=& \sum_{b=0}^1 \sum_{j \in 2\mathbb{Z} + b, \; n \in \mathbb{Z}/N} 
c_b^{(r, s) }( 4n - j^2) e^{2\pi i n \tau + 2\pi i j z} .\\ \nonumber
& & \qquad\qquad \qquad 0 \leq r, s\leq N-1.
\end{eqnarray}
In particular if $N=2,3,5,7$ we can write the twisted elliptic genus of $K3$ as follows \cite{DJS,DJS1}:
\begin{eqnarray}\label{pAtwist}
F^{(0,0)}( \tau, z) & =& \frac{8}{N} A(\tau, z), \\ \nonumber
F^{(0, s)} (\tau, z) &=& \frac{8}{N(N+1) } A(\tau, z)  - \frac{2}{N+1} B(\tau, z) {\tilde E}_N(\tau) 
, \quad \hbox{for} \, 1\leq s\leq (N-1),  \\ \nonumber
F^{(r, rk)}  &=& \frac{8}{N(N+1) } A(\tau, z)  + 
\frac{2}{N(N+1)} {\tilde E}_N\left( \frac{\tau +k}{N} \right) B(\tau, z),  \\ \nonumber
& &  \hbox{for}\, 1\leq r\leq (N-1), 1\leq k\leq (N-1)
\end{eqnarray}
where $B(\tau,z)=\frac{\theta_1^2(\tau,z)}{\eta^6(\tau)}$. $B(\tau,z)$ is a weak Jacobi form of weight $-2$ and index 1.

This quantity twisted elliptic genus is an ingredient in computing the statistical entropy of black holes in $\mathcal{N}=4$ type IIA/B string theory in 4 space-time dimensions \cite{DJS,DVV}.

Another instance where these modular functions $\tilde E_N$ is observed is in the context of computing new supersymmetric index in heterotic strings \cite{HM,CCL,SS}. 
 compactified on order $N$ orbifolds of $K3\times T^2$ and $E_8\times E_8$. 
This is defined by 
\begin{equation}
{\mathcal Z}_{\rm new} 
(q, \bar q )  = \frac{1}{\eta^2(\tau)  } {\rm Tr}_{R} 
\left(    F e^{i \pi  F } q^{L_0 - \frac{c}{24} } \bar q^{\bar L_0 -\frac{ \bar c}{24}}   \right)  \ . 
\end{equation}
Here the  trace  
is performed  over the Ramond sector in  the internal CFT with central 
charges $(c, \bar c) = ( 22, 9 ) $.  
$F$ refers to  the  world sheet fermion number of the right moving 
${\mathcal N}=2$ supersymmetric internal CFT. 
For a general embedding in case of an order 2 orbifold of $K3$ the new super-symmetric index can be given by,
\begin{eqnarray}\label{newindnon}
\mathcal{Z}_{\rm new}&=&-\frac{1}{\eta^{24}}\left\{ 2 \Gamma^{(0,0)}_{2, 2} 
E_4E_6  \right. \\ \nonumber
& & +\Gamma^{(0,1)}_{2,2} 
\left[(E_6+2{\tilde E}_2(\tau) E_4 ) \left(\hat{b}{\tilde E}_2^2(\tau)+ (\frac{2}{3} - \hat b)  E_4\right)\right]\\ \nn
&&+\Gamma^{(1,0)}_{2,2}
\left [\left(E_6-{\tilde E}_2(\frac{\tau}{2}) E_4\right)\left({\frac{\hat b}{4}}{\tilde E}_2^2(\frac{\tau}{2} )+ 
( \frac{2}{3} - \hat b) E_4\right) \right]\\ \nn
&&\left. +\Gamma^{(1,1)}_{2,2}
\left[ \left(E_6-{\tilde E}_2(\frac{\tau+1}{2})E_4\right) 
\left({\frac{\hat b}{4}}{\tilde E}_2^2(\frac{\tau+1}{2})+ ( \frac{2}{3} - \hat b) E_4\right) \right] \right\} .
\end{eqnarray}
where $\Gamma^{(r,s)}_{2,2}$ is a lattice sum depending on the moduli of the torus $T^2$ in the heterotic theory and $\hat b$ is a rational number.
Similarly for other orbifolds of $K3$ one can still write the $\mathcal{Z}_{\rm new}$ in terms of Eisenstein series of $SL(2,\mathbb{Z})$ and $\Gamma_0(N)$ \cite{ACJD}. This index is used to calculate the gauge and gravitational coupling corrections in heterotic strings which predicts the existence of different Calabi Yau 3-folds using heterotic-type II string duality \cite{ACJD1}.

\section{{\bf Notations and Preliminaries}}
For $k\ge 2$, even we denote the space of all modular forms on $SL(2,\mathbb{Z})$ by $M_k$ and the space of all cusp forms on $SL(2,\mathbb{Z})$ by $S_k$ and  the space of all modular forms on $\Gamma_0(N)$ by $M_k(\Gamma_0(N))$ and the space of all cusp forms on $\Gamma_0(N)$ by $S_k(\Gamma_0(N))$. 

 \noindent

Let $f(z)$ be a holomorphic function in the upper half plane and $\gamma=\left(\begin{matrix}
a & b\\c & d
\end{matrix}\right)\in GL_2(\mathbb{Q})$ then 
\[f|_k \gamma=\det(\gamma)^{k/2}(cz+d)^{-k}f(\frac{az+b}{cz+d}).\]

\noindent
With $\omega_N:=\left(\begin{matrix}0 & -1 \\ N & 0\end{matrix}\right) $ denoting the usual Fricke involution on $M_k$, one can readily verify the fact that $\tilde{E}_N|_2\; \omega_N=-\tilde{ E}_N$.

\noindent
Let $F_\Gamma$ be the fundamental domain (see \cite{TMA} for more details and definition) of $\Gamma:=SL_2(\mathbb{{Z}})$ given by the region in the upper half plane such as $\tau=x+iy$ where $-1/2\leq x<1/2$ and $x^2+y^2\geq 1$,
then $F_N$, the fundamental domain for $\Gamma_0(N)$  for a prime $N$ can be given as \cite{TMA}
\be\label{fundom}
F_N=F_{\Gamma}\,\cup\,{ {\bigcup}}_{k=0}^{N-1} ST^k (F_{\Gamma}),
\ee
where $\{\{ST^s\}_{s=0}^{N-1}\cup {\rm Id}\}$ is a set of coset representatives of $\Gamma_0(N)$ in $\Gamma$.

\noindent
{\bf The valence formula:} Let $f \neq 0$ be a modular form of weight $k\ge 2$, even in $SL(2,\mathbb{Z})$, $v_p(f)$ is the order of the zero of $f$ at point $p\in F_{\Gamma}$ and $\rho=e^{2\pi i/3}$ then the valence formula gives,
\be
v_{\infty}(f)+\frac{1}{2}v_i(f)+\frac{1}{3}v_{\rho}(f)+\sum_{p\in \Gamma, p\neq i, \rho}v_p=\frac{k}{12}.
\ee
From the above valence formula one obtains that in the fundamental domain of the full modular group the only zeros of $E_4$ and $E_6$ lie at $\tau=e^{2\pi i/3}$ and $\tau=i$ respectively.
\begin{equation}
E_6(i)=0,\qquad E_4(e^{2\pi i/3})=0.
\end{equation}
\section{{\bf Proof of Theorem \ref{theorem}}}\label{proof}
We shall begin this section with the following lemmas regarding the zeros of Eisenstein series which are essential to prove the \thmref{theorem}.

\begin{lem}\label{lem1}
	If $\tau\in F_N$ and
	$E_6(\tau)=0$ then either $\tau=i$ or $\tau=\frac{-1}{i+s}$. Similarly if $\tau\in F_N$ and $E_4(\tau)=0$ then either $\tau=e^{2\pi i/3}$ or $\tau=\frac{-1}{e^{2\pi i/3}+s}$ where $0\le s\le N-1$ and these are the only possible zeros of $E_4$ and $E_6$ in $F_N$.
	\end{lem}
 \proof The lemma follows from the expression of fundamental domain $F_N$ (\ref{fundom}) and the zeros of $E_4$ and $E_6$ in $F_\Gamma$.

\begin{lem}\label{poly}
	There exists a polynomial expression for $N=2,3,5,7$
	\[(\tilde{ E}_N(\tau))^{N+1}=\sum_{i=0 }^{N-1}a_{i} (\tilde{E}_N(\tau))^{i} m_{N+1-i}(\tau),\]
	where $a_{i}$ are constants and $m_{N+1-i} \in M_{2(N+1-i)}(SL_2(\mathbb{Z}))$.
	These polynomials are given by,
	\begin{eqnarray} \label{e2}
	\tilde{ E}_2(\tau)^3 &=& \frac{1}{4}E_6(\tau)+\frac{3}{4}E_4 \tilde{ E}_2(\tau), \\ \label{e3}
	\tilde{E}_3(\tau)^4 &=& \frac{1}{27}E_4^2(\tau)+\frac{8}{27}E_6(\tau) \tilde{E}_3(\tau)+\frac{2}{3}\tilde{E}_3(\tau)^2 E_4(\tau), \\ \nn
	\tilde{E}_5(\tau)^6 &=& \frac{1}{3125}E_6^2(\tau)+\frac{24}{3125}E_6 E_4 \tilde{E}_5(\tau)+\frac{9}{125} E_4 ^2 \tilde{E}_5(\tau)^2+\frac{8}{25} E_6\\ \label{e5} &&\tilde{E}_5(\tau)^3+\frac{3}{5}\tilde{E}_5(\tau)^4 E_4(\tau), \\ \nn
	\tilde{E}_7(\tau)^8 &=& \frac{1}{7^7}E_4^4+ \frac{48}{7^7}E_4^2 E_6 \tilde{E}_7+\left(\frac{64}{64827} E_6^2+\frac{92}{453789} E_4^3 \right) \tilde{E}_7^2+\frac{32}{2401}E_4 E_6 \tilde{E}_7^3\\ \label{e7}
	&&+\frac{30}{2401}E_4^2\tilde{E}_7^4+ \frac{16}{49}E_6 \tilde{E}_7^5+\frac{4}{7}\tilde{E}_7^6 E_4.
	\end{eqnarray}
\end{lem}

\proof Considering the $q$ expansion of both sides of the above equations upto Strum's bound the lemma follows.
The program  zeros.nb is given with this paper \cite{AC} for reference.

\proof[{\bf Proof of Theorem \ref{theorem} :}]
\noindent
From \lemref{lem1} we know the location of all zeros of $E_4$ and $E_6$ in $F_N$. Let us denote these sets of zeros in $F_N$ as $\mathbb{L}_{4,N}$ and $\mathbb{L}_{6,N}$ respectively. 
Also let us denote the set of all zeros of $\tilde{E}_N$ in $F_N$ as $\tilde{\mathbb{L}}_N$. From \lemref{poly} we observe that
\begin{equation*}
\tilde{\mathbb{L}}_2 \subseteq \mathbb{L}_{6,2}, \qquad \tilde{\mathbb{L}}_3 \subseteq \mathbb{L}_{4,3}, \qquad \tilde{\mathbb{L}}_5 \subseteq \mathbb{L}_{6,5}, \qquad \tilde{\mathbb{L}}_7 \subseteq \mathbb{L}_{4,7}.
\end{equation*}
Let us chose an element $\omega_{4,N}$ from $\mathbb{L}_{4,N}$ for $N=3,7$ and $\omega_{6,N}$ from $\mathbb{L}_{6,N}$ for $N=2,5$.
Now we re-write the resulting polynomial equations (\ref{e2})$-$(\ref{e7}) as a product of two factors as follows:
\begin{eqnarray} \label{ee2}
\tilde{ E}_2(\omega_{6,2})\left(\tilde{ E}_2(\omega_{6,2})^2-\frac{3}{4}E_4(\omega_{6,2})\right)  &=& 0, \\ \label{ee3}
\tilde{E}_3(\omega_{4,3})\left(\tilde{E}_3(\omega_{4,3})^3 -\frac{8}{27}E_6(\omega_{4,3})\right)  &=& 0, \\ \label{ee5}
\tilde{E}_5(\omega_{6,5})^2 \left(\tilde{E}_5(\omega_{6,5})^4 -\frac{3}{5}\tilde{E}_5(\omega_{6,5})^2 E_4(\omega_{6,5}) -\frac{9}{125} E_4(\omega_{6,5}) ^2\right) &=& 0, \\ \label{ee7}
\tilde{E}_7^2(\omega_{4,7})\left(\tilde{E}_7(\omega_{4,7})^6  - \frac{16}{49}E_6(\omega_{4,7}) \tilde{E}_7(\omega_{4,7})^3-\frac{64}{64827} E_6(\omega_{4,7})^2 \right) &=& 0.
\end{eqnarray}

    \noindent
     Since the right hand side of the above set of equations are zero so at least one of these factors must be zero.
 Now we need to check the numerical values of these factors to determine the location of zeros of $\tilde{ E}_N$ in $F_N$.

   So we estimate the bounds of $\tilde{E}_N(\omega)$, where $\omega=\omega_{6,N}$ for $N=2,5$ and $\omega=\omega_{4,N}$ for $N=3,7$. We write the $q$ expansion of $\tilde{ E}_N$ as follows:
  \begin{equation*}
  \tilde{ E}_N(\tau)=\sum_{n=0}^{m}a_n q^n+\sum_{n=m+1}^{\infty}a_n q^n=\sum_{n=0}^{m}a_n q^n+R(m,q).
  \end{equation*}
 Now we have,
  \begin{eqnarray} \label{bn}
 |R(m,q)|= |\sum_{n=m+1}^{\infty}a_n (x+iy)^n| \le  \sum_{n=m+1}^{\infty}|a_n| (|x|+|y|)^n \\ \nn
 \le \sum_{n=m+1}^{\infty}b_n (|x|+|y|)^n < \sum_{n=m+1}^{\infty}b_n (|x|+|y|+\epsilon)^n
  \end{eqnarray}
  where, $b_n=(N+1)n(n+1)/2> (N+1) \sigma_1(n)> |a_n|,\; q=x+iy, \; \epsilon>0$. {\footnote{$\epsilon$ is needed to approximate $|x|$ and $|y|$ in the computer program. Also $\sum_{n=1}^{m}a_n q^n$ is done up to machine precision.}}
  
  \noindent
  Note that our choice of $b_n$ is such that $\sum_{n=m+1}^{\infty}b_n (|x|+|y|+\epsilon)^n$ can be estimated exactly in terms of $|x|+|y|+\epsilon$ and $m$. Also note that  $\sum_{n=m+1}^{\infty}b_n (|x|+|y|+\epsilon)^n$ is convergent only if $|x|+|y|+\epsilon<1$.
  
  \noindent
  In the region $|x|+|y|+\epsilon<1$ we estimate the numerical value of $\tilde E_N(\omega)$ up to first $m+1$ terms in the $q$ expansion and using equation (\ref{bn}) we estimate the upper-bound of $|R(m,q)|$ in $F_N$. From this we have estimated an upper bound of $\tilde{ E}_N(\omega)$ for the values of $\omega$ as $\omega=\tau$ as in table \ref{table}) and lower bound of $|\tilde{ E}_N(\omega)|$ for the rest of the values of $\omega\in \mathbb{L}_{4,N}$ for $N=3,7$ and $\omega\in \mathbb{L}_{6,N}$ for $N=2,5$ (see table \ref{num} in the appendix for reference).
   
   \noindent
However if $|x|+|y|\ge 1$ (or close to 1) \footnote{This happens at $\omega\in \mathbb{L}_{6,5} $ at $\omega=\frac{-1}{i+4}$ and $\omega\in \mathbb{L}_{4,7}$ at $\omega=\frac{-1}{e^{2\pi i/3}+5}$, $\frac{-1}{e^{2\pi i/3}+6}$. At these points we can evaluate the result at $\omega=\frac{-1}{i-1},\frac{-1}{e^{2\pi i/3}-2}, \frac{-1}{e^{2\pi i/3}-1}$ respectively. These results can be related back to the fundamental domain points using equation (\ref{fun}).} we use the following results
  \begin{eqnarray}
  \tilde{ E}_N(\frac{\tau+s}{N})&=&\tilde{ E}_N(\frac{\tau+s-N}{N}),\\ \nn
  \tilde{ E}_N(-\frac{1}{\tau+s})&=&-\frac{(\tau+s)^2}{N}\tilde{ E}_N(\frac{\tau+s}{N}),
   \end{eqnarray}
 which is true for any $k\in \mathbb{Z}$. This implies that
 \be\label{fun}
  \tilde{ E}_N(-\frac{1}{\tau+s-N})=\frac{(\tau+s-N)^2}{(\tau+s)^2} \tilde{ E}_N(-\frac{1}{\tau+s}),
  \ee
  where $0\le s\le N-1$.
  Putting $\omega=-\frac{1}{\tau+s}$ in (\ref{fun}) we have,
  \be\label{funw}
  \tilde{ E}_N(\frac{1}{1/\omega+N})=(1/\omega+N)^2\omega^2 \tilde{ E}_N(\omega).
  \ee
  \noindent
  Note that if $\tau\in F_\Gamma$, then for $0\le s\le N-1$, we have $\frac{-1}{\tau+s}\in F_N$ and $\frac{-1}{\tau+s-N}\notin F_N$. So if $\omega=\omega_{4,N}$ for $N=3,7$ and $\omega=\omega_{6,N}$ for $N=2,5$ then the equation (\ref{funw}) implies that there is a relation between the value of $\tilde{ E}_N$ at a point outside $\mathbb{L}_{4,N}$ (respectively $\mathbb{L}_{6,N}$) and a point inside $\mathbb{L}_{4,N}$ (respectively $\mathbb{L}_{6,N}$). So at the point $\tau=\frac{1}{1/\omega+N}$ where $|\Re(e^{2\pi i\tau})|+|\Im(e^{2\pi i\tau})|<1$ (and not very close to 1) we repeat the process as before. Thus we can estimate the bounds of $|\tilde{ E}_N(\tau)|$ which in turn puts the bounds on $|\tilde{ E}_N(\omega)|$ (numerical estimates can be found in table \ref{num}).
 
 \noindent
 {\bf Estimates:}
 Now we shall list the estimates of $|\tilde{ E}_N(\omega)|$ upto $q^{m}$ terms in the expansion of $\tilde{ E}_N(\omega)$ where $\omega \in \mathbb{L}_{4,N}$ for $N=3,7$ and $\omega \in \mathbb{L}_{6,N}$ for $N=2,5$. 
 \begin{table}[H]
 	\renewcommand{\arraystretch}{-0.5}
 	\begin{center}
 		\vspace{0.5cm}
 		\begin{tabular}{|c|c|c|c|}
 			\hline
 			& & & \\
 			$N$ & $\omega$ & $|\tilde{E}_N(\omega)|$ & $R^>(100,q)$ \\
 			& & & \\
 			\hline
 			& & & \\
 			& & & \\
 			& & & \\
 			2 & $i$ & 1.0449 & $10^{-269}$  \\
 			& $\frac{-1}{i+1}$ & $0.0$ & $10^{-132}$ \\
 			\hline
 			& & &  \\
 			& & & \\
 			& & & \\
 			3 & $e^{2\pi i/3}$ & 0.948674 & $10^{-226}$ \\
 			& $\frac{-1}{e^{2\pi i/3}+2}$  & $3.14185\times 10^{-77} $ &$  10^{-73}$\\
 			\hline
 			& & & \\
 			& & & \\
 			& & & \\
 			5 & $i$ &1.01127 & $10^{-269}$  \\
 			& $\frac{-1}{i+1}$ & 0.77254 &  \\
 			& $\frac{-1}{i+2}$ &$2.96211\times 10^{-53}$ & $10^{-33}$ \\
 			& $\frac{-1}{i+3}$ & $1.65345\times 10^{-25}$ & $10^{-12}$ \\
 			& $\frac{-1}{i-1}$ & 0.77254 & $10^{-269}$\\
 			\hline
 			& & & \\
 			& & & \\
 			& & & \\
 			7 & $e^{2\pi i/3}$ &0.98289 & $10^{-226}$ \\
 			& $\frac{-1}{e^{2\pi i/3}+2}$ &0.615423 & $10^{-73}$   \\
 			& $\frac{-1}{e^{2\pi i/3}+3}$ & $2.63251\times 10^{-32}$ & $10^{-16}$ \\
 			& $\frac{-1}{e^{2\pi i/3}+4}$ &2.66684 & $10^{-7}$ \\
 			& $\frac{-1}{e^{2\pi i/3}-2}$ &$-2.63251\times 10^{-32}$ & $10^{-12}$ \\
 			& $\frac{-1}{e^{2\pi i/3}-1}$ & 0.615423  & $10^{-73}$ \\
 			\hline
 		\end{tabular}
 	\end{center}
 	\vspace{-0.2cm}
 	\caption{Numerical estimates of $|\tilde{ E}_N(\omega)|$ upto $q^{100}$ terms in the expansion of $\tilde{ E}_N(\omega)$ using machine precision up to 1000 digits in Mathematica and $R^>(100,q)$ is an upper bound on the remainder term $R(100,q)$}\label{num}
 	\renewcommand{\arraystretch}{-0.5}
 \end{table}
 
 For $m=100$ in all the above cases $|R(m,q)|<10^{-7}$. Using these bounds we get an upper bound of $|\tilde{ E}_N(\omega)|$ for $\omega=\tau$ as in table \ref{table} which can be given by $10^{-7}<10^{-2}$ so the second factors in equations (\ref{ee2}) to (\ref{ee7}) can not be zero at these points, so $\tilde{ E}_N(\omega)=0$. We also obtain the lower bound of $|\tilde{ E}_N(\omega)|$ for rest of the points in $\mathbb{L}_{4,N}$ for $N=3,7$ and $\mathbb{L}_{6,N}$ for $N=2,5$ and the bound is $>10^{-2}$, this shows that $\tilde{ E}_N\ne 0$ at those points. So the only points in $\tilde{\mathbb{L}}_N$ are the ones mentioned in table \ref{table}. For detailed calculations used above see the  mathematica file zeros.nb attached with the paper.

\noindent
Thus we obtain the points where $\tilde{ E}_N$ can be zero. 
   Now to prove that these are the only zeros of $\tilde{ E}_N$ in $F_N$ as given in table \ref{table} we must check that at these points the second factor will not become zero. We shall do this as follows:
 
 \noindent
  From the definition of $E_4$ and $E_6$ [(\ref{e4}) and (\ref{e6})] we see that $E_4(i)>1$ and $E_6(e^{2\pi i/3})>1$. Hence, from the modular transformation properties of $E_4$ and $E_6$ $|E_4(\omega_{6,2})|>1,\;\; |E_4(\omega_{6,5})|>1 $ and $|E_6(\omega_{4,2})|>1, \;\; |E_6(\omega_{4,7})|>1 $.
  So one can easily see that the second factor in equations (\ref{ee2}) to (\ref{ee7}) can only become zero if $|\tilde{ E}_N(\tau)|>10^{-2}$. So when the second factor would be zero that would never give a zero of $\tilde{ E}_N(\omega)$ and vice versa. Now from the positive bounds of $|\tilde{ E}_N|$ (see the table \ref{num}) we find all the zeros of $\tilde{ E}_N(\omega)$ which are given as:
  \bea \label{zeroen}
  \tilde{E}_2 \left(\frac{-1}{i+1}\right) &=&0 ,\\ \nn
  \tilde{E}_3 \left(\frac{-1}{e^{2\pi i/3}+2}\right) &=& 0, \\ \nn
  \tilde{E}_5 \left(\frac{-1}{i+2} \right) =\quad  \tilde{E}_5 \left(\frac{-1}{i+3} \right)&=&0, \\ \nn
  \tilde{E}_7 \left(\frac{-1}{e^{2\pi i/3}+3} \right)=\quad \tilde{E}_7 \left(\frac{-1}{e^{2\pi i/3}+5} \right)  &=&0.
  \eea
  
So these are the only possible zeros of $\tilde{ E}_N$ in $F_N$.
This completes the proof of \thmref{theorem} .

\section*{{\bf Remarks}}
\begin{enumerate}

\item The cusps of $\Gamma_0(N)$ are 0 and $i\infty$ for prime $N$. However from the $q$ expansions it is obvious that $\tilde{E}_N(i\infty)=1$. Also using Fricke involutions one can see that
\[\lim_{\epsilon\rightarrow 0}\tilde{E}_N(\frac{-1}{i\epsilon})=\lim_{\epsilon\rightarrow 0} \frac{(-(i\epsilon)^2)}{N}\tilde{E}_N(\frac{i\epsilon}{N}),\]
 so $\tilde{E}_N$ is non-zero at the cusps.
So the only possible zeros of $\tilde{ E}_N$ are at equations (\ref{zeroen}) in the fundamental domain $F_N$.

\item This method does not easily generalize to $\tilde{ E}_{11}$ (or higher prime numbers) because of the presence of cusp forms of weight 2 in $\Gamma_0(11)$. Also it doesn't easily generalize to composite numbers.
However using moonshine symmetry the twisted elliptic genus of $K3$ is known to exist for all 26 conjugacy classes of $M_{24}$ where higher $N$ Eisenstein series are present \cite{GHV,EH}.
\item One may try to generalize the method to different types of modular functions where there is a possibility of finding a polynomial relation in terms of modular functions whose zeros are known.

\item The equation (\ref{e2}) was used in finding the twisted elliptic genus, new supersymmetric index and gauge threshold corrections at one loop for the non-standard embedding of heterotic string compactified on $K3\times T^2$ where the $K3$ was orbifolded with a $\mathbb{Z}_2$ automorphism corresponding to the 2A conjugacy class of Mathieu group $M_{24}$ and a $1/2$ shift on one of the circles of $T^2$ \cite{ACJD}.

\end{enumerate}

\section{{\bf Acknowledgments}}
	The author thanks Council of Scientific and Industrial Research (CSIR) for the funding. We thank the anonymous referee for suggestions and comments. This work would not have been possible without the help and guidance of Prof. Justin R. David (Department of Centre for High Energy Physics, Indian Institute of Science) and Prof. Soumya Das (Department of Mathematics, Indian Institute of Science). The author also thanks Ritwik Pal and Pramath Anamby (Department of Mathematics, Indian Institute of Science) for useful discussions.


\begin{thebibliography}{9}


\bibitem[BS]{bs}
A. El Basraoui, and A. Sebbar,
\emph{Zeros of the Eisenstein series $E_2$}
. Proc. Amer. Math. Soc. 138 (2010),
no. 7, 2289–2299.

 \bibitem[CCL]{CCL}
G.~Lopes~Cardoso, G.~Curio, and D.~Lust, {\it {Perturbative couplings and
		modular forms in N=2 string models with a Wilson line}},  {\em Nucl.Phys.}
{\bf B491} (1997) 147--183, [\href{http://arxiv.org/abs/hep-th/9608154}{{\tt
		hep-th/9608154}}].

\bibitem[CD]{ACJD}
A. Chattopadhyaya and  J. R. David, \emph{ ${N}=2$ heterotic string compactifications on orbifolds of $K3\times T^2$},
JHEP 1701 (2017) 037.

\bibitem[CD1]{ACJD1}
A.Chattopadhyaya and J.R. David, {\emph {Gravitational couplings in ${\mathcal N}=2$ string
		compactifications and Mathieu Moonshine},  {JHEP 1805 (2018) 211}}
\url{https://arxiv.org/abs/1712.08791}.

\bibitem[DJS]{DJS}
J.~R. David, D.~P. Jatkar, and A.~Sen, {\it {Product representation of Dyon
		partition function in CHL models}},  {\em JHEP} {\bf 06} (2006) 064,
[\href{http://arxiv.org/abs/hep-th/0602254}{{\tt hep-th/0602254}}].

	\bibitem[DJS1]{DJS1}
J.~R. David and A.~Sen, {\it {CHL Dyons and Statistical Entropy Function from
		D1-D5 System}},  {\em JHEP} {\bf 11} (2006) 072,
[\href{http://arxiv.org/abs/hep-th/0605210}{{\tt hep-th/0605210}}].

\bibitem[DVV]{DVV}
R.~Dijkgraaf, E.~P. Verlinde, and H.~L. Verlinde, {\it {Counting dyons in N=4
		string theory}},  {\em Nucl. Phys.} {\bf B484} (1997) 543--561,
[\href{http://arxiv.org/abs/hep-th/9607026}{{\tt hep-th/9607026}}].

\bibitem[EH]{EH}
T.~Eguchi and K.~Hikami, {\it {Note on twisted elliptic genus of $K3$
		surface}},  {\em Phys. Lett.} {\bf B694} (2011) 446--455,
[\href{http://arxiv.org/abs/1008.4924}{{\tt arXiv:1008.4924}}].

\bibitem[EOTY]{EOTY}
T.~Eguchi, H.~Ooguri, A.~Taormina, and S.-K. Yang, {\it {Superconformal
		Algebras and String Compactification on Manifolds with SU(N) Holonomy}},
{\em Nucl. Phys.} {\bf B315} (1989) 193.

\bibitem[GHV]{GHV}
M.~R. Gaberdiel, S.~Hohenegger, and R.~Volpato, {\it {Mathieu twining
		characters for K3}},  {\em JHEP} {\bf 09} (2010) 058,
[\href{http://arxiv.org/abs/1006.0221}{{\tt arXiv:1006.0221}}].


	\bibitem[GS]{gar1}
S. Garthwaite, L. Long, H. Swisher, and S. Treneer, \emph{Zeros of classical Eisenstein series and recent developments}.

\bibitem[GSLT]{gslt}
  S. Garthwaite, L. Long, H. Swisher, and S. Treneer, \emph{Zeros of some level 2 Eisenstein series},  Proc. Amer. Math. Soc.,138(2010), 467-480.

 \bibitem[HM]{HM}
 J.~A. Harvey and G.~W. Moore, \emph{{Algebras, BPS states, and strings}},
 \href{https://doi.org/10.1016/0550-3213(95)00605-2}{\emph{Nucl. Phys.}
 	{\bfseries B463} (1996) 315--368},
 [\href{https://arxiv.org/abs/hep-th/9510182}{{\ttfamily hep-th/9510182}}].
 
 
 

	
	\bibitem[RSD]{RSD} F. K. C. Rankin and H. P. F. Swinnerton-Dyer, \emph{On the zeros of Eisenstein series
	} Bull. London Math. Soc.2
	(1970), 169–170
	
	\bibitem[RW]{woody}
	R. Wood, and M. P. Young. \emph{Zeros of the weight two Eisenstein Series}. Journal of Number Theory 143 (2014): 320-333.
	
	\bibitem[SJ]{shi07}
	J. Shigezumi, \emph{On the zeros of the Eisenstein series for
		$\Gamma_0 ^*(5)$ and $\Gamma_0 ^*(7)$}, Kyushu
	J. Math. 61
	(2007), no. 2, 527–549.
	
	\bibitem[SS]{SS}
	S.~Stieberger, \emph{{(0,2) heterotic gauge couplings and their M theory
			origin}}, \href{https://doi.org/10.1016/S0550-3213(98)00770-6}{\emph{Nucl.
			Phys.} {\bfseries B541} (1999) 109--144},
	[\href{https://arxiv.org/abs/hep-th/9807124}{{\ttfamily hep-th/9807124}}].
	

	\bibitem[TMA]{TMA}
	T. M. Apostol \emph{Modular Functions and Dirichlet Series in Number Theory}, Springer New York, 1976.
	
	

\end{thebibliography}
\end{document}